\documentclass[12pt]{article}
\usepackage{amssymb}
\usepackage{amsfonts}
\usepackage{times}
\usepackage{mathptmx}
\usepackage{amsmath}
\usepackage[usenames]{color}
\usepackage{mathrsfs}
\usepackage{amsfonts}
\usepackage{amssymb,amsmath}
\usepackage{CJK}
\usepackage{cite}
\usepackage{cases}
\usepackage{amsthm}
\def\cl{\centerline}

\def\vs{\vspace*}

\def\Z{\mathbb{Z}}

\def\C{\mathbb{C}}

\def\QED{\hfill$\Box$}
\def\ni{\noindent}

\numberwithin{equation}{section}
\newtheorem{theo}{Theorem}[section]
\newtheorem{defi}[theo]{Definition}

\newtheorem{lemm}[theo]{Lemma}
\newtheorem{remark}[theo]{Remark}
\newtheorem{examp}[theo]{Example}
\newtheorem{prop}[theo]{Proposition}

\begin{document}
\begin{center}
{\bf\large Generalized conformal derivations of Lie conformal algebras}
\footnote {Supported by the National Natural Science Foundation of China (No.~11431010, 11371278, 11501515) and the Zhejiang Provincial Natural Science Foundation of China (No.~LQ16A010011).

$^{\,\S}$Corresponding author: Y.~Hong.
}
\end{center}

\cl{Guangzhe Fan$^{\,*,\,\dag}$, Yanyong Hong$^{\,\S}$, Yucai Su$^{\,*,\,\ddag}$
}

\cl{\small $^{\,*}$Department of Mathematics, Tongji University, Shanghai 200092,
P. R. China}
\cl{\small $^{\,\S}$Department of Science, Zhejiang Agriculture and Forestry University,}
\cl{\small Hangzhou 311300, Zhejiang, P. R. China}

\cl{\small $^{\,\dag}$yzfanguangzhe@126.com}
\cl{\small $^{\,\S}$hongyanyong2008@yahoo.com}
\cl{\small $^{\,\ddag}$ycsu@tongji.edu.cn}
\vs{8pt}

{\small\footnotesize
\parskip .005 truein
\baselineskip 3pt \lineskip 3pt
\noindent{{\bf Abstract:} Let $R$ be a Lie conformal algebra. The purpose of this paper is to investigate the conformal derivation algebra $CDer(R)$, the conformal quasiderivation algebra $QDer(R)$ and the generalized conformal derivation algebra $GDer(R)$. The generalized conformal derivation algebra is a natural generalization of the conformal derivation algebra. Obviously, we have the following tower $CDer(R)\subseteq QDer(R)\subseteq GDer(R)\subseteq gc(R)$, where $gc(R)$ is the general Lie conformal algebra. Furthermore, we mainly research the connection of these generalized conformal derivations. Finally, the conformal $(\alpha,\beta,\gamma)$-derivations of Lie conformal algebras are studied. Moreover, we obtain some connections between several specific generalized conformal derivations and the conformal $(\alpha,\beta,\gamma)$-derivations. In addition, all conformal $(\alpha,\beta,\gamma)$-derivations of finite simple Lie conformal algebras are characterized.

\vs{5pt}


}}
\parskip .001 truein\baselineskip 6pt \lineskip 6pt

\section{Introduction}
Throughout this paper, we denote by $\C$ the set of complex numbers and $\C^\ast$ the set of non-zero complex numbers. We also assume that $\lambda,\mu$ and $\nu$ are indeterminate variables. For convenience, $\mathcal{A}$ will denote the ring $\C[\partial]$ of polynomials in the indeterminate $\partial$.
Moreover, if $V$ is a vector space, the space of polynomials of $\lambda$ with coefficients in $V$ is denoted by $V[\lambda]$.

It is well known that derivations and generalized derivations are very important subjects in the research of both Lie algebras and their generalizations.
All kinds of generalized derivations were studied in \cite{B,K3,LL1,LW,Z}. In recent years, many authors have investigated generalized derivations of nonassociative algebras, Lie algebras, Lie superalgebras and Lie color algebras in \cite{CML,H,LL,ZZ}. In \cite{LL}, G.~Leger and E.~Lucks investigated the generalized derivations of Lie algebras and their subalgebras. Moreover, they obtained many nice properties of the generalized derivation algebras and their subalgebras. Later on, many authors got many useful results similar to Lie algebras, such as Lie superalgebras, Lie color algebras.

Lie conformal algebra as a natural generalization of Lie algebra, was initially introduced by V. G. Kac as a formal language describing the singular part of the operator product expansion in conformal field theory. It is very useful for studying infinite dimensional Lie algebras satisfying the local property. Since then, the structure theories and representation theories of some Lie conformal algebras have been extensively undertaken in \cite{BKV,DK,FSW,FK,K1}. In the last few years, the conformal derivations of some concrete Lie conformal algebras have been studied by many authors in \cite{FSW,WCY}. We believe that the generalized conformal derivations of Lie conformal algebras would help to promote the development of structure theories of Lie conformal algebras. This is our motivation to present this paper.

In the present work, we would like to study the generalized conformal derivations of Lie conformal algebras. 
We mainly discuss the conformal derivation algebra $CDer(R)$, the center conformal derivation algebra $ZDer(R)$, the conformal quasiderivation algebra $QDer(R)$ and the generalized conformal derivation algebra $GDer(R)$ of a Lie conformal algebra $R$. It may be useful for investigating the structures of Lie conformal algebras. This is just the beginning of the study on generalized conformal derivations of Lie conformal algebras. Furthermore, we get many nice properties for Lie conformal algebras.

In \cite{NH,ZZ1}, ($\alpha$, $\beta$, $\gamma$)-derivations of Lie algebras and Lie superalgebras were studied. Moreover, they also obtained some properties of ($\alpha$, $\beta$, $\gamma$)-derivations of Lie (super)algebras. In particular, the spaces of ($\alpha$, $\beta$, $\gamma$)-derivations for simple Lie (super)algebras were presented. In this article, we will investigate the conformal ($\alpha$, $\beta$, $\gamma$)-derivations of Lie conformal algebras. And, the connection between generalized conformal derivations and the conformal $(\alpha,\beta,\gamma)$-derivations is studied. Moreover, conformal ($\alpha$, $\beta$, $\gamma$)-derivations of all finite simple Lie conformal algebras are characterized.

Let us briefly describe the structure of the article. In Section $2$, we review the basic notations which are used in this paper. In Section $3$, similar to the definitions of the generalized derivations of Lie algebras, some definitions of the generalized conformal derivations are presented first. Then, we study the connections of these generalized conformal derivations. Moreover, the generalized conformal derivations of Virasoro Lie conformal algebra are determined. In Section 4, we introduce the concept of conformal $(\alpha,\beta,\gamma)$-derivations and obtain some connections between the generalized conformal derivations and the conformal $(\alpha,\beta,\gamma)$-derivations. Finally, we classify the conformal $(\alpha,\beta,\gamma)$-derivations and characterize the conformal $(\alpha,\beta,\gamma)$-derivations of all finite simple Lie conformal algebras.

\section{Preliminaries}
In this section, for the reader's convenience, we shall summarize some basic facts about Lie conformal algebras used in this paper, see \cite{DK,K1,K2}.

\begin{defi}\label{250}\rm
A \emph{Lie conformal algebra} $R$ is an $\mathcal{A}$-module endowed with a $\C$-bilinear map
\begin{equation*}\label{25100}
R\otimes R\rightarrow R[\lambda], ~~~~~~a\otimes b\rightarrow [a{}\, _\lambda \, b],
\end{equation*}
satisfying the following three axioms~($a,b,c\in R$):
\begin{equation}\label{251}
\aligned
&Conformal~~sesquilinearity:~~~~[\partial a\,{}_\lambda \,b]=-\lambda[a\,{}_\lambda\, b],\ \ \ \
[a\,{}_\lambda \,\partial b]=(\partial+\lambda)[a\,{}_\lambda\, b];\\
&Skew~~symmetry:~~~~[a\, {}_\lambda\, b]=-[b\,{}_{-\lambda-\partial}\,a];\\
&Jacobi~~identity:~~~~[a\,{}_\lambda\,[b\,{}_\mu\, c]]=[[a\,{}_\lambda\, b]\,{}_{\lambda+\mu}\, c]+[b\,{}_\mu\,[a\,{}_\lambda \,c]].
\endaligned
\end{equation}
\end{defi}

The notions of a homomorphism, ideal and subalgebra of a Lie conformal
algebra are defined as usual. A Lie conformal algebra is
called \emph{finite} if it is finitely generated as an
$\mathcal{A}$-module. The \emph{rank} of a Lie conformal
algebra $R$ is its rank as an $\mathcal{A}$-module. A Lie conformal algebra $R$ is \emph{simple} if its only ideals are trivial and it is not abelian.

Suppose $A$ and $B$ are subspaces of a Lie conformal algebra $R$. We define $[A,B]$ as the $\mathbb{C}$-linear span of all $\lambda$-coefficients in the products $[a_\lambda b]$, where $a\in A$, $b\in B$. Moreover, if $B$ is an $\mathcal{A}$-submodule of $R$, by conformal sesquilinearity,
$[A,B]$ is also an $\mathcal{A}$-submodule of $R$.

There are two important examples of Lie conformal algebras.
\begin{examp}
The Virasoro Lie conformal algebra $\text{Vir}$ is the simplest nontrivial
example of Lie conformal algebras. It is defined by
$$\text{Vir}=\mathbb{C}[\partial]L, ~~[L_\lambda L]=(\partial+2\lambda)L.$$
\end{examp}
\begin{examp}

Let $\mathfrak{g}$ be a Lie algebra. The current Lie conformal
algebra associated to $\mathfrak{g}$ is defined by:
\begin{equation*}\label{25q2}
\text{Cur} \mathfrak{g}=\mathbb{C}[\partial]\otimes_{\mathbb{C}} \mathfrak{g}, ~~[a_\lambda b]=[a,b], \ \ \mbox{ \ for all $a,b\in \mathfrak{g}$.}
\end{equation*}
\end{examp}

It is shown in \cite{DK} that any finite simple Lie conformal algebra is either isomorphic to $\text{Vir}$ or $\text{Cur} \mathfrak{g}$ where $\mathfrak{g}$ is
a finite-dimensional simple Lie algebra.

\begin{defi}\rm\label{h1}
Suppose $R$ and $S$ are two $\mathcal{A}$-modules. A \emph{conformal linear map }from $R$ to $S$ is a $\C$-linear map $\phi_\lambda:R\longrightarrow S[\lambda]$ such that
\begin{equation*}\label{253}
\phi_\lambda(\partial a)=(\partial+\lambda)\phi_\lambda(a),\ \ \mbox{ \ for all $a\in R$.}
\end{equation*}
\end{defi}
\begin{remark}
Obviously, by Definition \ref{h1}, $\phi_\lambda$ does not depend on the choice of
the indeterminate variable $\lambda$. Therefore, for convenience, in the sequel, we also write this conformal linear map as $\phi_\mu$ for some other indeterminate variable $\mu$.
\end{remark}
Set $\phi_\lambda=\sum_{i=0}^n\lambda^i\phi_i$ where $\phi_i: R\rightarrow S$ is a $\C$-linear map for any $i\in\{0,\cdots,n\}$. Denote $ker (\phi)=\bigcap_{i=0}^n ker(\phi_i)$ and $Im (\phi)=\bigcup_{i=0}^n Im(\phi_i)$. Obviously,
$a\in ker (\phi)$ if and only if $\phi_\lambda(a)=0$ and $b\in Im(\phi)$ if
and only if it appears as a coefficient in a polynomial $\phi_\lambda(a)$ for
some $a\in R$.

Denote by $Chom(R,S)$ the space of conformal linear maps between $\mathcal{A}$-modules $R$ and $S$. For simplicity, we write $Cend(R)$ for $Chom(R,R)$. It can be made into an $\mathcal{A}$-module via $(\partial\phi)_\lambda(a)=-\lambda\phi_\lambda(a)$ for all $a\in R$.

\begin{defi}\rm
The $\lambda$-bracket on $Cend(R)$ given by
\begin{equation*}\label{255}
[\phi\,{}_\lambda \,\varphi]_{\mu}a=\phi_{\lambda}(\varphi_{\mu-\lambda}a)-\varphi_{\mu-\lambda}(\phi_{\lambda}a),\ \ \mbox{ \ for any $a\in R$,}
\end{equation*}
defines a Lie conformal algebra structure on $Cend(R)$. This is called the\emph{ general Lie conformal algebra }on $R$ and we denote it by $gc(R)$.
\end{defi}

\begin{defi}\rm
Let $R$ be a Lie conformal algebra. A conformal linear map $D_\lambda:R\longrightarrow R[\lambda]$ is called a \emph{conformal derivation} if
\begin{equation*}\label{254}
D_\lambda([a\,{}_\mu \,b])=[(D_\lambda a)\,{}_{\lambda+\mu} \,b]+[a\,{}_\mu \,(D_\lambda b)],\ \ \mbox{ \ for all $a,b\in R$.}
\end{equation*}
\end{defi}
It is easy to see that for any $a\in R$, the map $({\rm ad}_a)_\lambda$, defined by $({\rm ad}_a)_\lambda b= [a\, {}_\lambda\, b]$ for any $b\in R$, is a conformal derivation of $R$. All conformal derivations of this kind are called \emph{inner conformal derivations}. Denote by $CDer(R)$
and $CInn(R)$ the vector spaces of all conformal derivations and inner conformal derivations of $R$, respectively.

Obviously, $CInn(R)$ and $CDer(R)$ are conformal subalgebras of $gc(R)$ and
we have the following tower $CInn(R)\subseteq CDer(R)\subseteq gc(R)$.

\begin{defi}\rm
The center of Lie conformal algebra $R$ is $Z(R)=\{a\in R\,~|~[a_{\lambda}b]=0,~for~all~ b\in R\}$.
\end{defi}



\section{Generalized conformal derivation algebras}
Firstly, we will give some definitions and notations of the generalized conformal derivations of Lie conformal algebras.

An element $D_\lambda\in gc(R)$ is said to be a \emph{conformal quasiderivation} of a Lie conformal algebra $R$ if there exists an element  $D^{'}_\lambda\in gc(R)$
such that
\begin{equation}\label{e300}
[(D_\lambda a)\,{}_{\lambda+\mu} \,b]+[a\,{}_\mu \,(D_\lambda b)]=D^{'}_\lambda([a\,{}_\mu \,b]),\ \ \mbox{ \ for all $a,b\in R$.}
\end{equation}

Furthermore, if there exist two elements $D^{'}_\lambda,D^{''}_\lambda\in gc(R)$ such that
\begin{equation}\label{e301}
[(D_\lambda a)\,{}_{\lambda+\mu} \,b]+[a\,{}_\mu \,(D^{'}_\lambda b)]=D^{''}_\lambda([a\,{}_\mu \,b]),\ \ \mbox{ \ for all $a,b\in R$,}
\end{equation}
we call $D_\lambda$ a \emph{generalized conformal derivation}.

Denote by $QDer(R)$ and $GDer(R)$ the sets of conformal quasiderivations and generalized conformal derivations of $R$. By setting $D_\lambda^{'}=D_\lambda$ in (\ref{e301}), we can get $QDer(R)\subseteq GDer(R)$.

An element $D_\lambda\in gc(R)$ is called a\emph{ conformal centroid }of $R$, if it satisfies the following equality
\begin{equation}\label{303}
D_\lambda([a\,{}_\mu \,b])=[(D_\lambda a)\,{}_{\lambda+\mu} \,b]=[a\,{}_\mu \,(D_\lambda b)],\ \ \mbox{ \ for all $a,b\in R$.}
\end{equation}

If $D_\lambda\in gc(R)$ satisfies the following condition:
\begin{equation}\label{304}
[(D_\lambda a)\,{}_{\lambda+\mu} \,b]=[a\,{}_\mu \,(D_\lambda b)],\ \ \mbox{ \ for all $a,b\in R$,}
\end{equation}
$D_\lambda$ is called a \emph{conformal quasicentroid} of $R$.
From now on, for convenience, we denote by $C(R)$ and $QC(R)$ the sets of conformal centroids and quasicentroids. It is easy to see that $C(R)\subseteq QC(R)$.

Denote by $ZDer(R)$ the set of   elements $D_\lambda\in gc(R)$ satisfying
\begin{equation}\label{305}
[(D_\lambda a)\,{}_{\mu} \,b]=D_\lambda([a\,{}_\mu \, b])=0,\ \ \mbox{ \ for all $a,b\in R$.}
\end{equation}
It is easy to verify that $ZDer(R)$ is a Lie conformal ideal of $GDer(R)$ (see Proposition \ref{3002}).
We also have $D_\lambda a \in Z(R)[\lambda]$. If $Z(R)=0$, then $ZDer(R)=0$.

For the centerless Lie conformal algebra $R$, we have the following tower:
\begin{equation}\label{306}
R\cong CInn(R)\subseteq CDer(R)\subseteq QDer(R)\subseteq GDer(R)\subseteq gc(R).
\end{equation}
Moreover, when $R$ is abelian, it is easy to see $ZDer(R)=C(R)=QC(R)=CDer(R)=QDer(R)=GDer(R)=gc(R)$.

It is known that $CDer(Vir)=CInn(Vir)$ (see \cite{DK}). Then, we present a characterization of $C(Vir)$, $QC(Vir)$, $CDer(Vir)$,
$QDer(Vir)$ and $GDer(Vir)$.
\begin{prop}\label{po1}
For the Virasoro Lie conformal algebra, we can get
\begin{itemize}\parskip-3pt
\item[\rm(1)]$GDer(Vir)=QDer(Vir)=\{D_\lambda=a_0(\lambda)+a_1(\lambda)\partial~|~a_0(\lambda),
a_1(\lambda)\in \mathbb{C}[\lambda]\}$.
\item[\rm(2)]$C(Vir)=QC(Vir)=0.$
\end{itemize}
\end{prop}
\ni\ni{\it Proof.}\ \
Let $D_\lambda\in GDer(Vir)$. Then, there exist two elements $D_\lambda^{'}$,
$D_\lambda^{''}\in gc(Vir)$ such that
\begin{equation}\label{eq1}
[(D_\lambda L)_{\lambda+\mu}L]+[L_\mu(D_\lambda^{'}L)]=D_\lambda^{''}([L_\mu L]).
\end{equation}
Assume that $D_\lambda(L)=\sum_{i=0}^m a_i(\lambda)\partial^iL$,
$D_\lambda^{'}(L)=\sum_{i=0}^nb_i(\lambda)\partial^iL$ and
$D_\lambda^{''}(L)=\sum_{i=0}^sc_i(\lambda)\partial^iL$, where $a_i(\lambda)$,
$b_i(\lambda)$, $c_i(\lambda)\in \mathbb{C}[\lambda]$ and $a_m(\lambda)\neq 0$,
$b_n(\lambda)\neq 0$, $c_s(\lambda)\neq 0$.
Taking them into (\ref{eq1}) and comparing the coefficients of $L$, we can get
\begin{equation}\label{eq2}
(\partial+2\lambda+2\mu)\sum_{i=0}^ma_i(\lambda)(-\lambda-\mu)^i
+(\partial+2\mu)\sum_{i=0}^nb_i(\lambda)(\mu+\partial)^i
=(\partial+\lambda+2\mu)\sum_{i=0}^sc_i(\lambda)\partial^i.
\end{equation}
Obviously, by comparing the highest degrees of $\partial$ in both two sides of (\ref{eq2}),
$s$ must be equal to $n$.
If $n>1$, letting $\mu=0$ in (\ref{eq2}) and comparing the coefficients of $\partial^{n+1}$ and $\partial^n$, we obtain
\begin{equation}\label{eq3}
b_n(\lambda)=c_n(\lambda),~~b_{n-1}(\lambda)=c_{n-1}(\lambda)+c_n(\lambda)\lambda.
\end{equation}
Then, according to (\ref{eq3}) and by comparing the coefficients of $\partial^n$
in (\ref{eq2}), we get $n\mu b_n(\lambda)=0$. This shows that $b_n(\lambda)=0$, a contradiction. Then, by comparing the coefficients of $\mu$, we obtain $m\leq 1$. Therefore, $D_\lambda(L)=(a_0(\lambda)
+a_1(\lambda)\partial)L$, $D_\lambda^{'}(L)=(b_0(\lambda)
+b_1(\lambda)\partial)L$ and $D_\lambda^{''}(L)=(c_0(\lambda)+c_1(\lambda)\partial)L$.
Plugging them into (\ref{eq2}), we can get $b_1(\lambda)=c_1(\lambda)=a_1(\lambda)$,
$b_0(\lambda)=a_0(\lambda)$ and $c_0(\lambda)=2(a_0(\lambda)-a_1(\lambda)\lambda)$.

Thus, we obtain
\begin{equation}
D_\lambda(L)=D_\lambda^{'}(L)=(a_0(\lambda)+a_1(\lambda)\partial)L,~~~
D_\lambda^{''}(L)=(2(a_0(\lambda)-a_1(\lambda)\lambda)+a_1(\lambda)\partial)L.
\end{equation}
Therefore, by the definitions of conformal quasiderivation and generalized conformal derivation, we get $GDer(Vir)=QDer(Vir)=\{D_\lambda=a_0(\lambda)+a_1(\lambda)\partial~|~a_0(\lambda),
a_1(\lambda)\in \mathbb{C}[\lambda]\}$.

Similarly, it is easy to see $C(Vir)=QC(Vir)=0$. \QED
\begin{remark}
By Proposition \ref{po1}, we can see that $CDer(Vir)\neq GDer(Vir)=QDer(Vir)$.
\end{remark}

Next, we begin to study the connections of these generalized conformal derivations.
\begin{prop}\label{3001}
$GDer(R)$, $QDer(R)$ and $C(R)$ are Lie conformal subalgebras of $gc(R)$.
\end{prop}
\ni\ni{\it Proof.}\ \ Here, we only prove that for any $D_{1\,{}_\lambda},D_{2\,{}_\lambda}\in GDer(R)$, $[D_1\,{}_\lambda \,D_2]\,{}_\mu \in GDer(R)[\lambda]$. The other two cases can be proved similarly.

According to the definition of $gc(R)$, we get
\begin{equation*}\label{31000}
[D_1\,{}_\lambda \,D_2]\,{}_\mu=D_{1\,{}_\lambda}D_{2\,{}_{\mu-\lambda}}-D_{2\,{}_{\mu-\lambda}}D_{1\,{}_\lambda}.
\end{equation*}
Since ${D_1}_{\lambda}$ and ${D_2}_{\lambda}\in GDer(R)$,
there exist four elements  ${D^{'}_1}_\lambda$, ${D^{'}_2}_\lambda$, ${D^{''}_1}_\lambda$, ${D^{''}_2}_\lambda\in gc(R)$ such that
\begin{eqnarray*}
&&[({D_1}_\lambda a)\,{}_{\lambda+\mu} \,b]+[a\,{}_\mu \,({D_1^{'}}_\lambda b)]={D_1^{''}}_\lambda([a\,{}_\mu \,b]),\\
&&[({D_2}_\lambda a)\,{}_{\lambda+\mu} \,b]+[a\,{}_\mu \,({D_2^{'}}_\lambda b)]={D_2^{''}}_\lambda([a\,{}_\mu \,b]), \text{ for all $a,b\in R$.}
\end{eqnarray*}
Then, we obtain
\begin{eqnarray*}\!\!\!\!\!\!\!\!\!\!\!\!&\!\!\!\!\!\!\!\!\!\!\!\!\!\!\!&
\ \ \ \ \ \ \ \ [(D_1\,{}_\lambda \,(D_2\,{}_{\mu-\lambda}a))_{\mu+\nu}b]
=D^{''}_1\,{}_\lambda([(D_2\,{}_{\mu-\lambda}a)\,{}_{\nu+\mu-\lambda}b])-[(D_2\,{}_{\mu-\lambda} a)\,{}_{\nu+\mu-\lambda}(D^{'}_1\,{}_\lambda b)]
\nonumber\\\!\!\!\!\!\!\!\!\!\!\!\!&\!\!\!\!\!\!\!\!\!\!\!\!\!\!\!&
\ \ \ \ \ \ \ \ \ \ \ \ \ \ \ \ \ \ \ \ \ \ \ \ \ \ \ \ \ \ \ \ \ \ \ \ \ \ \ \ \ \ \ \, =D^{''}_1\,{}_\lambda(D^{''}_2\,{}_{\mu-\lambda}([a\,{}_{\nu}b])-[a\,{}_{\nu}(D^{'}_2\,{}_{\mu-\lambda}b)])
\nonumber\\\!\!\!\!\!\!\!\!\!\!\!\!&\!\!\!\!\!\!\!\!\!\!\!\!\!\!\!&
\ \ \ \ \ \ \ \ \ \ \ \ \ \ \ \ \ \ \ \ \ \ \ \ \ \ \ \ \ \ \ \ \ \ \ \ \ \ \ \ \ \ \ \ \ \ \, -D^{''}_2\,{}_{\mu-\lambda}([a\,{}_{\nu}(D^{'}_1\,{}_\lambda b)])+[a\,{}_{\nu}(D^{'}_2\,{}_{\mu-\lambda}(D^{'}_1\,{}_\lambda b))].
\end{eqnarray*}
Similarly, we have
\begin{eqnarray*}\!\!\!\!\!\!\!\!\!\!\!\!&\!\!\!\!\!\!\!\!\!\!\!\!\!\!\!&
\ \ \ \ \ \ \ \ [(D_{2\,{}_{\mu-\lambda}}(D_{1\,{}_\lambda}a))_{\mu+\nu}b]
 =D^{''}_2\,{}_{\mu-\lambda}(D^{''}_1\,{}_{\lambda}([a\,{}_{\nu}b])-[a\,{}_{\nu}(D^{'}_1\,{}_{\lambda}b)])
\nonumber\\\!\!\!\!\!\!\!\!\!\!\!\!&\!\!\!\!\!\!\!\!\!\!\!\!\!\!\!&
\ \ \ \ \ \ \ \ \ \ \ \ \ \ \ \ \ \ \ \ \ \ \ \ \ \ \ \ \ \ \ \ \ \ \ \ \ \ \ \ \ \ \ \ \ \ \, -D^{''}_1\,{}_{\lambda}([a\,{}_{\nu}(D^{'}_2\,{}_{\mu-\lambda} b)])+[a\,{}_{\nu}(D^{'}_1\,{}_{\lambda}(D^{'}_2\,{}_{\mu-\lambda} b))].
\end{eqnarray*}
Thus, it follows at once that
\begin{eqnarray*}\!\!\!\!\!\!\!\!\!\!\!\!&\!\!\!\!\!\!\!\!\!\!\!\!\!\!\!&
\ \ \ \ \ \ \ \ [([D_1\,{}_\lambda \,D_2]\,{}_{\mu}a)_{\mu+\nu}b]+[a_{\nu}([D^{'}_1\,{}_\lambda \,D^{'}_2]\,{}_{\mu}b)]
=[D^{''}_1\,{}_\lambda D^{''}_2]\,{}_{\mu}([a\,{}_{\nu}b]).
\end{eqnarray*}
Therefore, $[D_1\,{}_\lambda \,D_2]\,{}_\mu \in GDer(R)[\lambda]$, i.e. $GDer(R)$ is a Lie conformal subalgebra of $gc(R)$.\QED

\begin{prop}\label{3002}
$ZDer(R)$ is a Lie conformal ideal of $GDer(R)$.
\end{prop}
\ni\ni{\it Proof.}\ \
Firstly, it is readily seen that $ZDer(R)\subseteq GDer(R)$.


Next, we only need to prove that $[D_1\,{}_\lambda \,D_2]\,{}_\mu\in ZDer(R)[\lambda]$ for any $D_1\,{}_\lambda\in ZDer(R)$, $D_2\,{}_\lambda\in GDer(R)$.
Since ${D_2}_{\lambda}\in GDer(R)$,
there exist two elements ${D^{'}_2}_\lambda$, ${D^{''}_2}_\lambda\in gc(R)$ such that
\begin{eqnarray*}
[({D_2}_\lambda a)\,{}_{\lambda+\mu} \,b]+[a\,{}_\mu \,({D_2^{'}}_\lambda b)]={D_2^{''}}_\lambda([a\,{}_\mu \,b]), \text{ for all $a,b\in R$.}
\end{eqnarray*}
Then, we obtain
\begin{eqnarray*}\!\!\!\!\!\!\!\!\!\!\!\!&\!\!\!\!\!\!\!\!\!\!\!\!\!\!\!&
\ \ \ \ \ \ \ \ [([D_1\,{}_\lambda \,D_2]\,{}_{\mu}a)_{\nu}b]
=[(D_1\,{}_\lambda \,(D_2\,{}_{\mu-\lambda}a))_{\nu}b]-[(D_2\,{}_{\mu-\lambda} \,(D_1\,{}_\lambda a))_{\nu}b]
\nonumber\\\!\!\!\!\!\!\!\!\!\!\!\!&\!\!\!\!\!\!\!\!\!\!\!\!\!\!\!&
\ \ \ \ \ \ \ \ \ \ \ \ \ \ \ \ \ \ \ \ \ \ \ \ \ \ \ \ \ \ \ \ \ \ \ \ \ \, =0-(D_2^{''}\,{}_{\mu-\lambda}([(D_1\,{}_{\lambda}a)_{\nu+\lambda-\mu}b])-[(D_1\,{}_{\lambda}a)\,{}_{\nu+\lambda-\mu}(D^{'}_2\,{}_{\mu-\lambda}b)])
\nonumber\\\!\!\!\!\!\!\!\!\!\!\!\!&\!\!\!\!\!\!\!\!\!\!\!\!\!\!\!&
\ \ \ \ \ \ \ \ \ \ \ \ \ \ \ \ \ \ \ \ \ \ \ \ \ \ \ \ \ \ \ \ \ \ \ \ \ \,
=0-(D^{''}_2\,{}_{\mu-\lambda}(0)-0)
\nonumber\\\!\!\!\!\!\!\!\!\!\!\!\!&\!\!\!\!\!\!\!\!\!\!\!\!\!\!\!&
\ \ \ \ \ \ \ \ \ \ \ \ \ \ \ \ \ \ \ \ \ \ \ \ \ \ \ \ \ \ \ \ \ \ \ \ \ \,
=0.
\end{eqnarray*}
And, it is easy to see that $[D_1\,{}_\lambda \,D_2]_\mu([a_\nu b])=0$ for any
$a$, $b\in R$.

Hence, we have completed the proof.\QED

\begin{remark}
Let $R$ be a Lie conformal algebra which is free and of finite rank as an
$\mathcal{A}$-module. It is known that $gc(R)$ is a simple Lie conformal algebra.
Therefore, if $GDer(R)=gc(R)$, by Proposition \ref{3002}, in this case, $ZDer(R)$ is equal to either $0$ or $gc(R)$.
\end{remark}

\begin{prop}\label{3003}
If $Z(R)=0$, then $C(R)$ is a commutative Lie conformal subalgebra of $GDer(R)$.
\end{prop}
\ni\ni{\it Proof.}\ \
Obviously, we have $C(R)\subseteq GDer(R)$ by their definitions.

Let ${D_1}_\lambda$, ${D_2}_\lambda\in C(R)$. By some computations, we deduce that
\begin{eqnarray*}\!\!\!\!\!\!\!\!\!\!\!\!&\!\!\!\!\!\!\!\!\!\!\!\!\!\!\!&
\ \ \ \ \ \ \ \ [(D_1\,{}_\lambda \,(D_2\,{}_{\mu-\lambda}a))\,{}_{\nu}b]
=D_1\,{}_\lambda([(D_2\,{}_{\mu-\lambda}a)\,{}_{\nu-\lambda}b])
\nonumber\\\!\!\!\!\!\!\!\!\!\!\!\!&\!\!\!\!\!\!\!\!\!\!\!\!\!\!\!&
\ \ \ \ \ \ \ \ \ \ \ \ \ \ \ \ \ \ \ \ \ \ \ \ \ \ \ \ \ \ \ \ \ \ \ \ \ \ \ \ \ \ \,=D_1\,{}_\lambda([a\,{}_{\nu-\mu}(D_2\,{}_{\mu-\lambda}b)])
\nonumber\\\!\!\!\!\!\!\!\!\!\!\!\!&\!\!\!\!\!\!\!\!\!\!\!\!\!\!\!&
\ \ \ \ \ \ \ \ \ \ \ \ \ \ \ \ \ \ \ \ \ \ \ \ \ \ \ \ \ \ \ \ \ \ \ \ \ \ \ \ \ \ \,=[(D_1\,{}_\lambda a)\,{}_{\nu+\lambda-\mu}(D_2\,{}_{\mu-\lambda}b)]
\nonumber\\\!\!\!\!\!\!\!\!\!\!\!\!&\!\!\!\!\!\!\!\!\!\!\!\!\!\!\!&
\ \ \ \ \ \ \ \ \ \ \ \ \ \ \ \ \ \ \ \ \ \ \ \ \ \ \ \ \ \ \ \ \ \ \ \ \ \ \ \ \ \ \,=[(D_2\,{}_{\mu-\lambda}\,(D_1\,{}_\lambda a))\,{}_{\nu}b].
\end{eqnarray*}
It follows that\begin{equation*}\label{31200}
[([D_1\,{}_\lambda\,D_2]\,{}_{\mu}a)\,{}_{\nu}b]=0.
\end{equation*}
Due to the arbitrary of $b$ and $Z(R)=0$, we obtain that
\begin{equation*}\label{3100}
[D_1\,{}_\lambda\,D_2]\,{}_{\mu}a=0.
\end{equation*}
This means $[D_1\,{}_\lambda\,D_2]\,{}_{\mu}=0$.

Now we get the expected result.\QED

\begin{lemm}\label{r3004}
The following properties hold:
\begin{itemize}\parskip-3pt
\item[\rm(1)]$C(R)+CDer(R)\subseteq QDer(R)$.
\item[\rm(2)]$C(R)\subseteq QDer(R)\cap QC(R)$.
\item[\rm(3)]$GDer(R)=QDer(R)+QC(R)$.
\end{itemize}
\end{lemm}
\ni\ni{\it Proof.}\ \
(1)~For any $D_\lambda\in C(R)$, it follows that
\begin{equation*}\label{3043}
[(D_\lambda a)\,{}_{\lambda+\mu} \,b]+[a\,{}_\mu \,(D_\lambda b)]=2D_\lambda([a\,{}_\mu \,b]).
\end{equation*}

Thus, we can choose $D^{'}_\lambda=2 D_\lambda$ such that  $D_\lambda\in QDer(R)$.

Due to the definition of $CDer(R)$ and $QDer(R)$, we have seen that $CDer(R)\subseteq QDer(R)$. Therefore, (1) holds.

(2)~It is obvious that $C(R)\subseteq QC(R)$. Furthermore, according to (1), then (2) holds.

(3)~Assume $D_\lambda\in GDer(R)$. Then, there exist two elements $D^{'}_\lambda,D^{''}_\lambda\in gc(R)$ such that
\begin{equation}\label{q301}
[(D_\lambda a)\,{}_{\lambda+\mu} \,b]+[a\,{}_\mu \,(D^{'}_\lambda b)]=D^{''}_\lambda([a\,{}_\mu \,b]),\ \ \mbox{ \ for all $a,b\in R$.}
\end{equation}
Set $D_1\,{}_\lambda=\frac{D_\lambda+D^{'}_\lambda}2$ and $D_2\,{}_\lambda=\frac{D_\lambda-D^{'}_\lambda}2$, we only need to show that $D_1\,{}_\lambda\in QDer(R)$ and $D_2\,{}_\lambda\in QC(R)$.

By skew symmetry and (\ref{q301}), for any $a$, $b\in R$, we get
\begin{equation}\label{w1}
[b_{-\lambda-\mu-\partial}(D_\lambda a)]+[(D_\lambda^{'}b)_{-\mu-\partial}a]=D_\lambda^{''}([b_{-\mu-\partial}a]).
\end{equation}
By the conformal sesquilinearity and replacing $-\lambda-\mu-\partial$ by $\mu^{'}$ in (\ref{w1}), we obtain
\begin{equation}\label{w2}
[b_{\mu^{'}}(D_\lambda a)]+[(D_\lambda^{'}b)_{\lambda+\mu^{'}}a]=D_\lambda^{''}([b_{\mu^{'}}a]).
\end{equation}
Then, by changing the place of $a$, $b$ and replacing $\mu^{'}$ by $\mu$ in (\ref{w2}), we have
\begin{equation}\label{w3}
[a_{\mu}(D_\lambda b)]+[(D_\lambda^{'} a)_{\lambda+\mu}b]=D_\lambda^{''}([a_{\mu}b]).
\end{equation}

For one thing, by (\ref{q301}) and (\ref{w3}), we have
\begin{equation}\label{3054}
[(D_1\,{}_\lambda a)\,{}_{\lambda+\mu} \,b]+[a\,{}_\mu \,(D_1\,{}_\lambda b)]=D^{''}_\lambda([a\,{}_\mu \,b]),\ \ \mbox{ \ for all $a,b\in R$.}
\end{equation}
Hence, $D_1\,{}_\lambda\in QDer(R)$.

For another,  by (\ref{q301}) and (\ref{w3}), we also obtain that for any $a$, $b\in R$,
\begin{equation}\label{3040}
[(D_2\,{}_\lambda a)\,{}_{\lambda+\mu} \,b]=[a\,{}_\mu \,(D_2\,{}_\lambda b)].
\end{equation}

Therefore, $D_2\,{}_\lambda\in QC(R)$.\QED

\begin{lemm}\label{r3005}
We have the following conclusions:
\begin{itemize}\parskip-3pt
\item[\rm(1)]$[CDer(R)\,{}_\lambda \,C(R)]_{\mu}\subseteq C(R)[\lambda]$.
\item[\rm(2)]$[QDer(R)\,{}_\lambda \,QC(R)]_{\mu}\subseteq QC(R)[\lambda]$.
\item[\rm(3)]$[QC(R)\,{}_\lambda \,QC(R)]_{\mu}\subseteq QDer(R)[\lambda]$.
\end{itemize}
\end{lemm}
\ni\ni{\it Proof.}\ \ (1)~Let $D_1\,{}_\lambda\in CDer(R)$, $D_2\,{}_\lambda\in C(R)$. Then we get
\begin{eqnarray*}\!\!\!\!\!\!\!\!\!\!\!\!&\!\!\!\!\!\!\!\!\!\!\!\!\!\!\!&
\ \ \ \ \ \ \ \ [(D_1\,{}_\lambda \,(D_2\,{}_{\mu-\lambda}a))_{\mu+\nu}b]
=D_1\,{}_\lambda([(D_2\,{}_{\mu-\lambda}a)\,{}_{\nu+\mu-\lambda}b])-[(D_2\,{}_{\mu-\lambda} a)\,{}_{\nu+\mu-\lambda}(D_1\,{}_\lambda b)]
\nonumber\\\!\!\!\!\!\!\!\!\!\!\!\!&\!\!\!\!\!\!\!\!\!\!\!\!\!\!\!&
\ \ \ \ \ \ \ \ \ \ \ \ \ \ \ \ \ \ \ \ \ \ \ \ \ \ \ \ \ \ \ \ \ \ \ \ \ \ \ \ \ \ \ \, =D_1\,{}_\lambda (D_2\,{}_{\mu-\lambda}([a\,{}_{\nu}b]))-D_2\,{}_{\mu-\lambda}([a\,{}_{\nu}(D_1\,{}_\lambda b)]).
\end{eqnarray*}
\begin{eqnarray*}\!\!\!\!\!\!\!\!\!\!\!\!&\!\!\!\!\!\!\!\!\!\!\!\!\!\!\!&
\ \ \ \ \ \ \ \ [(D_2\,{}_{\mu-\lambda} \,(D_1\,{}_\lambda a))_{\mu+\nu}b]
=D_2\,{}_{\mu-\lambda}([(D_1\,{}_\lambda a)\,{}_{\nu+\lambda}b])
\nonumber\\\!\!\!\!\!\!\!\!\!\!\!\!&\!\!\!\!\!\!\!\!\!\!\!\!\!\!\!&
\ \ \ \ \ \ \ \ \ \ \ \ \ \ \ \ \ \ \ \ \ \ \ \ \ \ \ \ \ \ \ \ \ \ \ \ \ \ \ \ \ \ \ \,
=D_2\,{}_{\mu-\lambda}(D_1\,{}_\lambda([a\,{}_{\nu}b]))-D_2\,{}_{\mu-\lambda}([a\,{}_{\nu}(D_1\,{}_\lambda b)]).
\end{eqnarray*}
Furthermore, one has
\begin{eqnarray*}\!\!\!\!\!\!\!\!\!\!\!\!&\!\!\!\!\!\!\!\!\!\!\!\!\!\!\!&
\ \ \ \ \ \ [([D_1\,{}_\lambda \,D_2]\,{}_{\mu}a)_{\mu+\nu}b]
=[((D_1\,{}_\lambda \,D_2\,{}_{\mu-\lambda}-D_2\,{}_{\mu-\lambda} \,D_1\,{}_\lambda)a)_{\mu+\nu}b]
\nonumber\\\!\!\!\!\!\!\!\!\!\!\!\!&\!\!\!\!\!\!\!\!\!\!\!\!\!\!\!&
\ \ \ \ \ \ \ \ \ \ \ \ \ \ \ \ \ \ \ \ \ \ \ \ \ \ \ \ \ \ \ \ \ \ \ \ \ \ \ \,
=(D_1\,{}_\lambda \,D_2\,{}_{\mu-\lambda}-D_2\,{}_{\mu-\lambda} \,D_1\,{}_\lambda)([a_{\nu}b])
\nonumber\\\!\!\!\!\!\!\!\!\!\!\!\!&\!\!\!\!\!\!\!\!\!\!\!\!\!\!\!&
\ \ \ \ \ \ \ \ \ \ \ \ \ \ \ \ \ \ \ \ \ \ \ \ \ \ \ \ \ \ \ \ \ \ \ \ \ \ \ \,
=[D_1\,{}_\lambda \,D_2]\,{}_{\mu}([a_{\nu}b]).
\end{eqnarray*}
Similarly, we can obtain
\begin{eqnarray*}\!\!\!\!\!\!\!\!\!\!\!\!&\!\!\!\!\!\!\!\!\!\!\!\!\!\!\!&
\ \ \ \ \ \ \ \ [([D_1\,{}_\lambda \,D_2]\,{}_{\mu}a)_{\mu+\nu}b]
=[a_{\nu}([D_1\,{}_\lambda \,D_2]\,{}_{\mu}b)].
\end{eqnarray*}

For (2), using the similar method of proving (1), we can conclude that (2) holds.

For (3), suppose that $D_1\,{}_\lambda, D_2\,{}_\lambda\in QC(R)$, we can obtain that
\begin{eqnarray*}\!\!\!\!\!\!\!\!\!\!\!\!&\!\!\!\!\!\!\!\!\!\!\!\!\!\!\!&
\ \ \ \ \ \ \ \ [([D_1\,{}_\lambda \,D_2]\,{}_{\mu}a)_{\mu+\nu}b]
=[(D_1\,{}_\lambda \,(D_2\,{}_{\mu-\lambda}a))_{\mu+\nu}b]-[(D_2\,{}_{\mu-\lambda} \,(D_1\,{}_\lambda a))_{\mu+\nu}b]
\nonumber\\\!\!\!\!\!\!\!\!\!\!\!\!&\!\!\!\!\!\!\!\!\!\!\!\!\!\!\!&
\ \ \ \ \ \ \ \ \ \ \ \ \ \ \ \ \ \ \ \ \ \ \ \ \ \ \ \ \ \ \ \ \ \ \ \ \ \ \ \ \ \,
=[(D_2\,{}_{\mu-\lambda}a)_{\mu+\nu-\lambda}D_1\,{}_\lambda b)]-[(D_1\,{}_\lambda a)_{\nu+\lambda}(D_2\,{}_{\mu-\lambda}b)]
\nonumber\\\!\!\!\!\!\!\!\!\!\!\!\!&\!\!\!\!\!\!\!\!\!\!\!\!\!\!\!&
\ \ \ \ \ \ \ \ \ \ \ \ \ \ \ \ \ \ \ \ \ \ \ \ \ \ \ \ \ \ \ \ \ \ \ \ \ \ \ \ \ \,
=[a_{\nu}(D_2\,{}_{\mu-\lambda} \,(D_1\,{}_\lambda b))]-[a_{\nu}(D_1\,{}_\lambda \,(D_2\,{}_{\mu-\lambda}b))]
\nonumber\\\!\!\!\!\!\!\!\!\!\!\!\!&\!\!\!\!\!\!\!\!\!\!\!\!\!\!\!&
\ \ \ \ \ \ \ \ \ \ \ \ \ \ \ \ \ \ \ \ \ \ \ \ \ \ \ \ \ \ \ \ \ \ \ \ \ \ \ \ \ \,
=-[a_{\nu}([D_1\,{}_\lambda \,D_2]\,{}_{\mu}b)].
\end{eqnarray*}

Hence, we get
\begin{equation*}\label{3030001}
[([D_1\,{}_\lambda \,D_2]\,{}_{\mu}a)_{\mu+\nu}b]+[a_{\nu}([D_1\,{}_\lambda \,D_2]\,{}_{\mu}b)]=0([a_{\nu}b]).
\end{equation*}
This completes the proof.\QED

\begin{prop}\label{3006}
$QC(R)$+$[QC(R),QC(R)]$ is a Lie conformal ideal of $GDer(R)$.
\end{prop}
\ni\ni{\it Proof.}\ \ In view of Lemma \ref{r3004}~(3) and Lemma \ref{r3005}~(2),(3), it suffices to prove the statement.
\QED



\begin{prop}\label{3008}
If $Z(R)=0$, then $QC(R)$ is a Lie conformal algebra if and only if $[QC(R)\,{}_\lambda \,QC(R)]_{\mu}$ $=0.$ In particular, if $Z(R)=0$, then $C(R)$ centralizes $QC(R)$.
\end{prop}
\ni\ni{\it Proof.}\ \ $\Leftarrow$ It is nothing to prove.

$\Rightarrow$ Assume that $D_1\,{}_\lambda, D_2\,{}_\lambda\in QC(R)$, then we have $[D_1\,{}_\lambda \,D_2]\,{}_{\mu} \in QC(R)[\lambda]$.
Furthermore, we obtain that
\begin{equation*}\label{3703}
[([D_1\,{}_\lambda \,D_2]\,{}_{\mu}a)_{\mu+\nu}b]=[a_{\nu}([D_1\,{}_\lambda \,D_2]\,{}_{\mu}b)].
\end{equation*}

According to Lemma \ref{r3005}~(3), we obtain
\begin{equation*}\label{3803}
[([D_1\,{}_\lambda \,D_2]\,{}_{\mu}a)_{\mu+\nu}b]+[a_{\nu}([D_1\,{}_\lambda \,D_2]\,{}_{\mu}b)]=0([a_{\nu}b]).
\end{equation*}
Finally, it shows that
\begin{equation*}\label{3903}
[([D_1\,{}_\lambda \,D_2]\,{}_{\mu}a)_{\mu+\nu}b]=0.
\end{equation*}
According to $Z(R)=0$, this forces $[D_1\,{}_\lambda \,D_2]\,{}_{\mu}=0$.

This proves the proposition.\QED


\begin{prop}\label{3011}
If $D_\lambda\in C(R)$, then $ker(D)$ and $Im(D)$ are ideals of $R$.
\end{prop}
\ni\ni{\it Proof.}\ \
Assume that $D_\lambda\in C(R)$, $b\in ker(D)$. By the definition of $ker(D)$,
we know $D_\lambda(b)=0$. Then, for any $a\in R$, by the definition of $C(R)$, we have
\begin{eqnarray*}\!\!\!\!\!\!\!\!\!\!\!\!&\!\!\!\!\!\!\!\!\!\!\!\!\!\!\!&
\ \ \ \ \ \ \ \ D_\lambda([a\,{}_\mu \,b])=[(D_\lambda a)\,{}_{\lambda+\mu} \,b]=[a\,{}_{\mu} \,(D_\lambda b)]=0.
\end{eqnarray*}
Therefore, $[a_\mu b]\in ker(D)[\mu]$.

Suppose that $c\in Im(D)$. Then, we assume that $c$ appears as a coefficient in $D_\lambda d$ for some $d\in R$.
Therefore, we get
\begin{eqnarray*}\!\!\!\!\!\!\!\!\!\!\!\!&\!\!\!\!\!\!\!\!\!\!\!\!\!\!\!&
\ \ \ \ \ \ \ \ [(D_\lambda d)\,{}_{\mu} \,b]=D_\lambda([d\,{}_{\mu-\lambda} \,b]).
\end{eqnarray*}
Therefore, $[c_\mu b]\in Im(D)[\mu]$.

Hence, this proves $ker(D)$ and $Im(D)$ are ideals of $R$.\QED

Let $R$ be a centerless Lie conformal algebra. Assume that $R=S\oplus T$, where $S,T$ are ideals of $R$. If $D_\lambda\in GDer(R)$, it is easy to see that $D_\lambda(S)\subseteq S[\lambda]$, $D_\lambda(T)\subseteq T[\lambda]$. Furthermore, we get the following results.
\begin{prop}\label{3010}
\begin{itemize}\parskip-3pt
\item[\rm(1)]$GDer(R)=GDer(S)\oplus GDer(T)$.
\item[\rm(2)]$QDer(R)=QDer(S)\oplus QDer(T)$.
\item[\rm(3)]$C(R)=C(S)\oplus C(T)$.
\item[\rm(4)]$QC(R)=QC(S)\oplus QC(T)$.
\end{itemize}
\end{prop}

Similar to the investigation of quasiderivations of  Lie algebras, we shall deduce that the conformal quasiderivations of some Lie conformal algebras can be embedded as conformal derivations in some larger Lie conformal algebras.

Let $R$ be a Lie conformal algebra and $t$ an indeterminate. We define $\widetilde{R}$, canonically associated to $R$, as $\widetilde{R}=R[t\C[t]\diagup (t^3)]=Rt+Rt^2$. $\widetilde{R}$ can be endowed with a natural $\mathcal{A}$-module structure as follows: $\partial(a\otimes t^{i})=(\partial a )\otimes t^i $ for
any $a\in R$ and $i=1$, $2$. It is readily seen that $\widetilde{R}$ is a Lie conformal algebra with the $\lambda$-bracket:
$$[a\otimes t^{i}\, {}_\lambda\, b\otimes t^{j}]=[a\, {}_\lambda\, b]\otimes t^{i+j},$$ where $a,b\in R$, $i,j\in \Z$.
For simplicity, we write $at$ and $at^2$ in place of $a\otimes t$ and $a\otimes t^2$ respectively.

In the following of this section, we present an assumption: $R$ has
a decomposition of $\mathcal{A}$-modules, i.e.
 \begin{equation}\label{wq1}
 R=[R,R]\oplus U,~~~~\text{where $U$ is an $\mathcal{A}$-submodule of $R$.}
  \end{equation}
  Then $\widetilde{R}=Rt+Rt^2=Rt+[R,R]t^2+Ut^2$.
Now we define a map $f:QDer(R)\longrightarrow gc(\widetilde{R})$ satisfying $f(D_\lambda)(at+bt^2+ut^2)=D_\lambda(a)t+D^{'}_\lambda(b)t^2$,
where $a\in R$, $b\in [R,R]$, $u\in U$, and $D^{'}_\lambda$ is an element in $gc(R)$ associated with $D_\lambda$ in (\ref{e300}). Obviously, by the definition, $f(D_\lambda)$ is a conformal linear map.

\begin{theo}\label{3100000}
Let $R$ be the Lie conformal algebra satisfying (\ref{wq1}). Then, we have the following conclusions:
\begin{itemize}\parskip-3pt
\item[\rm(1)]$f$ is injective.
\item[\rm(2)]$f(D_\lambda)$ does not depend on the choice of $D^{'}_\lambda$.
\item[\rm(3)]$f(QDer(R))\subseteq CDer(\widetilde{R})$.
\end{itemize}
\end{theo}
\ni\ni{\it Proof.}\ \ (1)~If $f(D_1\,{}_\lambda)=f(D_2\,{}_\lambda)$, then we get
$f(D_1\,{}_\lambda)(at+bt^2+ut^2)=f(D_2\,{}_\lambda)(at+bt^2+ut^2)$.
This shows that
$D_1\,{}_\lambda(a)t+D^{'}_1\,{}_\lambda(b)t^2=D_2\,{}_\lambda(a)t+D^{'}_2\,{}_\lambda(b)t^2$.
Therefore, we conclude that $D_1\,{}_\lambda(a)=D_2\,{}_\lambda(a)$ for any $a\in R$.
Hence $D_1\,{}_\lambda=D_2\,{}_\lambda$, and $f$ is injective.

(2)~Assume that there exists $D^{''}_\lambda\in gc(R)$ such that $f(D\,{}_\lambda)(at+bt^2+ut^2)=D_\lambda(a)t+D^{''}_\lambda(b)t^2$.
Furthermore, for any $c$, $d\in R$, we have
\begin{equation}\label{310001}
[(D_\lambda c)\,{}_{\lambda+\mu} \,d]+[c\,{}_\mu \,(D_\lambda d)]=D^{'}_\lambda([c\,{}_\mu \,d]),
\end{equation}

\begin{equation}\label{310002}
[(D_\lambda c)\,{}_{\lambda+\mu} \,d]+[c\,{}_\mu \,(D_\lambda d)]=D^{''}_\lambda([c\,{}_\mu \,d]).
\end{equation}
It follows that $D^{'}_\lambda([c\,{}_\mu \,d])=D^{''}_\lambda([c\,{}_\mu \,d])$ for any $c$, $d\in R$.
Thus, for any $b\in [R,R]$, we get $D^{'}_\lambda(b)=D^{''}_\lambda(b)$.
Hence
\begin{eqnarray*}\!\!\!\!\!\!\!\!\!\!\!\!&\!\!\!\!\!\!\!\!\!\!\!\!\!\!\!&
\ \ \ \ \ \ \ \ f(D\,{}_\lambda)(at+bt^2+ut^2)=D_\lambda(a)t+D^{'}_\lambda(b)t^2
\nonumber\\\!\!\!\!\!\!\!\!\!\!\!\!&\!\!\!\!\!\!\!\!\!\!\!\!&
\ \ \ \ \ \ \ \ \ \ \ \ \ \ \ \ \ \ \ \ \ \ \ \ \ \ \ \ \ \ \ \ \ \ \ \ \ \ \ \ \ \ \ \ \,
=D_\lambda(a)t+D^{''}_\lambda(b)t^2,
\end{eqnarray*}
which shows that $f(D\,{}_\lambda)$ does not depend on the choice of $D^{'}_\lambda$.

(3)~According to the definition of $\widetilde{R}$, we obtain $[at^{i}\,{}_\mu \,bt^{j}]=0$ for $i+j\geq 3$.
Therefore, we only need to prove the following equality $f(D\,{}_\lambda)([at\,{}_\mu \,bt])=[(f(D\,{}_\lambda)at)\,{}_{\lambda+\mu} \,bt]+[at\,{}_\mu \,(f(D\,{}_\lambda) bt)]$.

By some computations, one has
\begin{eqnarray*}\!\!\!\!\!\!\!\!\!\!\!\!&\!\!\!\!\!\!\!\!\!\!\!\!\!\!\!&
\ \ \ \ \ \ \ \ f(D\,{}_\lambda)([at\,{}_\mu \,bt])=f(D\,{}_\lambda)([a\,{}_\mu \,b]t^{2})
\nonumber\\\!\!\!\!\!\!\!\!\!\!\!\!&\!\!\!\!\!\!\!\!\!\!\!\!&
\ \ \ \ \ \ \ \ \ \ \ \ \ \ \ \ \ \ \ \ \ \ \ \ \ \ \ \ \ \ \ \ \ \ \,
=D^{'}_\lambda([a\,{}_\mu \,b])t^{2}
\nonumber\\\!\!\!\!\!\!\!\!\!\!\!\!&\!\!\!\!\!\!\!\!\!\!\!\!&
\ \ \ \ \ \ \ \ \ \ \ \ \ \ \ \ \ \ \ \ \ \ \ \ \ \ \ \ \ \ \ \ \ \ \,
=([(D_\lambda a)\,{}_{\lambda+\mu} \,b]+[a\,{}_\mu \,(D_\lambda b)])t^{2}
\nonumber\\\!\!\!\!\!\!\!\!\!\!\!\!&\!\!\!\!\!\!\!\!\!\!\!\!&
\ \ \ \ \ \ \ \ \ \ \ \ \ \ \ \ \ \ \ \ \ \ \ \ \ \ \ \ \ \ \ \ \ \ \,
=[(f(D\,{}_\lambda)at)\,{}_{\lambda+\mu} \,bt]+[at\,{}_\mu \,(f(D\,{}_\lambda) bt)].
\end{eqnarray*}

Finally, we have completed the proof.\QED
\begin{remark}
Obviously, not all finite Lie conformal algebras satisfy the assumption.
For example, $R=\mathbb{C}[\partial]L\oplus \mathbb{C}[\partial]W$ is
a Lie conformal algebra of rank 2 with the $\lambda$-bracket as follows:
\begin{equation}
[L_\lambda L]=(\partial+2\lambda)L,~~[L_\lambda W]=\partial W,~~[W_\lambda W]=0.
\end{equation}
Then, $[R,R]=\mathbb{C}[\partial]L\oplus \mathbb{C}[\partial]\partial W$.
Therefore, there does not exist any $\mathcal{A}$-submodule $U$ such that
$R=[R,R]\oplus U$.

Of course, there are also many Lie conformal algebras satisfying the assumption.
All finite simple, semisimple Lie conformal algebras and $Cur \mathfrak{g}$ where
$\mathfrak{g}$ is a finite-dimensional Lie algebra are all such Lie conformal algebras. In addition, there are some solvable Lie conformal algebras satisfying the
assumption. For example, $R=\mathbb{C}[\partial]a\oplus \mathbb{C}[\partial]b$ is
a Lie conformal algebra of rank 2 with the $\lambda$-bracket as follows:
\begin{equation}
[a_\lambda a]=0,~~[a_\lambda b]=d(\lambda) b,~~[b_\lambda b]=0,
\end{equation}
where $d(\lambda)\in \mathbb{C}[\lambda]$. Then, if $d(\lambda)=0$,
$[R,R]=0$, i.e. $U=R$. If $d(\lambda)\neq 0$, $[R,R]=\mathbb{C}[\partial]b$.
Then, $U=\mathbb{C}[\partial]a$.
\end{remark}

\begin{theo}\label{31001}
Let $R$ be the Lie conformal algebra satisfying (\ref{wq1}). If $Z(R)=0$, and $\widetilde{R}$, $f$ are defined as above, then $CDer(\widetilde{R})=f(QDer(R))\oplus ZDer(\widetilde{R})$.
\end{theo}
\ni\ni{\it Proof.}\ \  From $Z(R)=0$, it follows at once that $Z(\widetilde{R})=Rt^2$. Suppose $\widetilde{D_\lambda}\in CDer(\widetilde{R})$, then we obtain $\widetilde{D_\lambda}(Z(\widetilde{R}))\subseteq Z(\widetilde{R})[\lambda]$. Thus, $\widetilde{D_\lambda}(Ut^2)\subseteq\widetilde{D_\lambda}(Z(\widetilde{R}))\subseteq Z(\widetilde{R})[\lambda]=Rt^2[\lambda]$.

Now define a conformal linear map:
$g_\lambda:Rt+[R,R]t^2+Ut^2\longrightarrow Rt^2[\lambda]$ by
\begin{equation}\label{111111}
g_\lambda(x)=\left\{\begin{array}{lll}
\widetilde{D_\lambda}(x)\cap Rt^2[\lambda]&\mbox{if \ }x\in Rt,\\
\widetilde{D_\lambda}(x)&\mbox{if \ }x\in Ut^2,\\
0&\mbox{if \ }x\in [R,R]t^2.
\end{array}\right.
\end{equation}

Obviously, we know that $g_\lambda([\widetilde{R},\widetilde{R}])\subseteq g_\lambda([R,R]t^2)=0$. Furthermore, we have $[g_\lambda(\widetilde{R})_\mu\widetilde{R}]\subseteq [Rt^2[\lambda]_\mu(Rt+Rt^2)]=0$.
Hence, we get $g_\lambda\in ZDer(\widetilde{R})$.

By the definition of $g_\lambda$, we have $(\widetilde{D_\lambda}-g_\lambda)(Rt)=\widetilde{D_\lambda}(Rt)\cap Rt[\lambda]\subseteq Rt[\lambda]$.
We also have $(\widetilde{D_\lambda}-g_\lambda)(Ut^2)=0$ and $(\widetilde{D_\lambda}-g_\lambda)([R,R]t^2)=\widetilde{D_\lambda}([\widetilde{R},\widetilde{R}])\subseteq [\widetilde{R},\widetilde{R}][\lambda]\subseteq [R,R]t^2[\lambda]$.

Therefore, there exist $D_\lambda,D'_\lambda\in gc(R)$ for all $a\in R$, $b\in [R, R]$ such that
$(\widetilde{D_\lambda}-g_\lambda)(at)=D_\lambda(a)t$,
$(\widetilde{D_\lambda}-g_\lambda)(bt^2)=D'_\lambda(b)t^2$.

Since $\widetilde{D_\lambda}-g_\lambda\in CDer(\widetilde{R})$ and by the definition of $CDer(\widetilde{R})$, we have $[(\widetilde{D_\lambda}-g_\lambda)(a_{1}t)\,{}_\mu \,a_{2}t]+[a_{1}t\,{}_\mu \,(\widetilde{D_\lambda}-g_\lambda)(a_{2}t)]=(\widetilde{D_\lambda}-g_\lambda)[a_{1}t\,{}_\mu \,a_{2}t]$. Thus, $[D_\lambda(a_{1})\,{}_\mu \,a_{2}]+[a_{1}\,{}_\mu \,D_\lambda(a_{2})]=D'_\lambda([a_{1}\,{}_\mu \,a_{2}])$. Therefore, $D_\lambda\in QDer(R)$.

Furthermore, we deduce $\widetilde{D_\lambda}\!-\!g_\lambda\!=\! f(D_\lambda)\!\subseteq\! f(QDer(R))$,
then $CDer(\widetilde{R})\!=\!f(QDer(R))\!+\!ZDer(\widetilde{R})$.

For any $\widetilde{D_\lambda}\in f(QDer(R))\cap ZDer(\widetilde{R})$, there exists an element $\widetilde{D_\lambda}\in QDer(R)$ such that $\widetilde{D_\lambda}=f(D_\lambda)$. Then, $\widetilde{D_\lambda}(at+bt+ut^2)=f(D_\lambda)(at+bt+ut^2)=D_\lambda(a)t+D'_\lambda(b)t^2$, for $a\in R$, $b\in [R,R]$. Since $\widetilde{D_\lambda}\in ZDer(\widetilde{R})$, then $\widetilde{D_\lambda}(at+bt+ut^2)\in Z(\widetilde{R})[\lambda]=Rt^2[\lambda]$.
Therefore, $\widetilde{D_\lambda}(a)=0$ for any $a\in R$. Finally, we conclude $\widetilde{D_\lambda}=0$.

Thus, we have completed the proof. \QED

\section{Conformal $(\alpha,\beta,\gamma)$-derivations}
Throughout this section, we are concerned with the conformal $(\alpha,\beta,\gamma)$-derivations of Lie conformal algebras. The conformal $(\alpha,\beta,\gamma)$-derivations are natural generalizations of conformal derivations.

Let $R$ be an arbitrary Lie conformal algebra. We call a conformal linear map $D_\lambda\in Cend R$ a \emph {conformal $(\alpha,\beta,\gamma)$-derivation} of $R$ if there exist $\alpha,\beta,\gamma\in \C$ such that for any $a,b\in R$, the following relation is satisfied:
\begin{equation}\label{4000}
\alpha D_\lambda([a\,{}_\mu \,b])=[(\beta D_\lambda a)\,{}_{\lambda+\mu} \,b]+[a\,{}_\mu \,(\gamma D_\lambda b)].
\end{equation}
For any given $\alpha,\beta,\gamma\in \C$, we denote the set of all conformal $(\alpha,\beta,\gamma)$-derivations by $CDer_{(\alpha,\beta,\gamma)}R$, i.e.

\begin{equation*}\label{40000}
CDer_{(\alpha,\beta,\gamma)}R=\{D_\lambda\in Cend R  \,~|~\, \alpha D_\lambda([a\,{}_\mu \,b])=[(\beta D_\lambda a)\,{}_{\lambda+\mu} \,b]+[a\,{}_\mu \,(\gamma D_\lambda b)], \forall~a,b\in R\}.
\end{equation*}
Obviously, we have the following conclusions.
\begin{prop}\label{40001}\begin{itemize}\parskip-3pt
\item[\rm(1)]$CDer_{(1,1,1)}R=CDer(R)$.
\item[\rm(2)]$CDer_{(0,1,-1)}R=QC(R)$.
\item[\rm(3)]$CDer_{(1,0,0)}R\cap CDer_{(0,1,0)}R=ZDer(R)$.
\item[\rm(4)]$CDer_{(\alpha,\beta,\gamma)}R=CDer_{(\triangle\alpha,\triangle\beta,\triangle\gamma)}R$,\ \ \mbox{ \ for any $\triangle\in \C^*$.}
\end{itemize}
\end{prop}

Furthermore, we have the following proposition. The proposition will be used later on.

\begin{prop}\label{40002}
$CDer_{(\alpha,\beta,\gamma)}R=CDer_{(0,\beta-\gamma,\gamma-\beta)}R\cap CDer_{(2\alpha,\beta+\gamma,\beta+\gamma)}R$, where $\alpha,\beta,\gamma\in \C$.
\end{prop}
\ni\ni{\it Proof.}\ \
For any $D_\lambda\in CDer_{(\alpha,\beta,\gamma)}R$ and $a,b\in R$, we get the following equality:
\begin{equation}\label{4001}
\alpha D_\lambda([a\,{}_\mu \,b])=[(\beta D_\lambda a)\,{}_{\lambda+\mu} \,b]+[a\,{}_\mu \,(\gamma D_\lambda b)].
\end{equation}
With a similar discussion as that in the proof of Lemma \ref{r3004} (3), we can also get
\begin{equation}\label{et3}
\alpha D_\lambda([a_{\mu}b])=[a_{\mu}(\beta D_\lambda b)]+[(\gamma D_\lambda a)_{\lambda+\mu}b].
\end{equation}
Then, by (\ref{4001}) and (\ref{et3}), we get
\begin{equation}\label{4003}
0=(\beta-\gamma)([( D_\lambda a)\,{}_{\lambda+\mu} \,b]-[a\,{}_\mu \,( D_\lambda b)]),
\end{equation}
and
\begin{equation}\label{4004}
2\alpha D_\lambda([a\,{}_\mu \,b])=(\beta+\gamma)([( D_\lambda a)\,{}_{\lambda+\mu} \,b]+[a\,{}_\mu \,( D_\lambda b)]).
\end{equation}

Therefore, $CDer_{(\alpha,\beta,\gamma)}R\subseteq CDer_{(0,\beta-\gamma,\gamma-\beta)}R\cap CDer_{(2\alpha,\beta+\gamma,\beta+\gamma)}R$.
Similarly, from (\ref{4003}) and (\ref{4004}), we can also conclude that (\ref{4001}) holds.

Now, this proposition is proved.\QED

In the following, we would like to find the connection of these three parameters. Furthermore, we formulate that the three original parameters are depended on only one form from the following theorem.

Using the similar method in \cite{NH}, we have the following theorem.
\begin{theo}\label{40003}
For any $\alpha,\beta,\gamma\in \C$, there exists $\delta\in\C$ such that the subspace $CDer_{(\alpha,\beta,\gamma)}R\subseteq Cend R$ has only the following four cases:~(1)~$CDer_{(\delta,0,0)}R$, (2)~$CDer_{(\delta,1,-1)}R$, (3)~$CDer_{(\delta,1,0)}R$, (4)~$CDer_{(\delta,1,1)}R$.
\end{theo}

From above discussions, we get the main result of this section.
\begin{theo}\label{r40003}
For any $\alpha,\beta,\gamma\in \C$, then $CDer_{(\alpha,\beta,\gamma)}R$ is equal to one of the following subspaces of $Cend R$:
\begin{itemize}\parskip-3pt
\item[\rm(i)]$CDer_{(0,0,0)}R=Cend~R$.
\item[\rm(ii)]$CDer_{(1,0,0)}R=\{D_\lambda\in Cend~R  \,~|~\,  D_\lambda([R,R])=0\}$.
\item[\rm(iii)]$CDer_{(0,1,-1)}R=QC(R)$.
\item[\rm(iv)]$CDer_{(\delta,1,-1)}R=CDer_{(0,1,-1)}R\cap CDer_{(1,0,0)}R$.
\item[\rm(v)]$CDer_{(\delta,1,1)}R$, $\delta\in\C$.
\item[\rm(vi)]$CDer_{(\delta,1,0)}R=CDer_{(0,1,-1)}R\cap CDer_{(2\delta,1,1)}R$.
\end{itemize}
\end{theo}

Finally, we give a characterization of all conformal $(\alpha,\beta,\gamma)$-derivations of finite simple Lie conformal algebras.
\begin{theo}
(1) For the Virasoro Lie conformal algebra and any $\alpha$, $\beta$, $\gamma\in \mathbb{C}$, then $CDer_{(\alpha,\beta,\gamma)}Vir$ is equal to one of the following subspaces of $Cend~Vir$:
\begin{itemize}\parskip-3pt
\item[\rm(i)]$CDer_{(0,0,0)}Vir=Cend~Vir,$
\item[\rm(ii)]$CDer_{(1,0,0)}Vir=CDer_{(0,1,-1)}Vir=CDer_{(\delta,1,-1)}Vir=CDer_{(\delta,1,0)}Vir=0,$
\item[\rm(iii)]$CDer_{(\delta,1,1)}Vir=0$, if $\delta\neq 1$ and $\delta\neq 2$; $CDer_{(1,1,1)}Vir=CInn(Vir)$ and $CDer_{(2,1,1)}Vir=\{D_\lambda~|~D_\lambda(L)=a(\lambda)L$, \text{where $a(\lambda)\in \mathbb{C}[\lambda]$}\}.
\end{itemize}
(2) For the current Lie conformal algebra $A=Cur \mathfrak{g}$ when $\mathfrak{g}$ is a finite-dimensional simple Lie algebra and any $\alpha$, $\beta$, $\gamma\in \mathbb{C}$, then $CDer_{(\alpha,\beta,\gamma)}A$ is equal to one of the following subspaces of $Cend~A$:
\begin{itemize}\parskip-3pt
\item[\rm(i)]$CDer_{(0,0,0)}A=Cend~A$,
\item[\rm(ii)]$CDer_{(1,0,0)}A=CDer_{(\delta,1,-1)}A=0$ if $\delta\neq 0$,
\item[\rm(iii)]$CDer_{(0,1,-1)}A=\{D_\lambda~|~D_\lambda(a)=f(\lambda)a$, \text{where $a\in\mathfrak{g}$ and $f(\lambda)\in \mathbb{C}[\lambda]$.}\}
\item[\rm(iv)]$CDer_{(\delta,1,1)}A=0$ if $\delta\neq 1$ and $\delta\neq 2$; $CDer_{(1,1,1)}A=\{D_\lambda~|~D_\lambda(a)=P(\partial)(\partial+\lambda)a+d_\lambda(a),$
    \text{ where  $a\in\mathfrak{g}$, $P(\partial)\in \mathbb{C}[\partial]$ and
    $d_\lambda\in CInn(A)$,}\} and $CDer_{(2,1,1)}A=\{D_\lambda~|~D_\lambda(a)=f(\lambda)a$,\\ \text{where $a\in\mathfrak{g}$ and $f(\lambda)\in \mathbb{C}[\lambda]$}\}.
\item[\rm(v)] $CDer_{(\delta,1,0)}A=0$ if $\delta\neq 1$; $CDer_{(1,1,0)}A=\{D_\lambda~|~D_\lambda(a)=f(\lambda)a$, \text{where $a\in\mathfrak{g}$, $f(\lambda)\in \mathbb{C}[\lambda]$.}\}
\end{itemize}
\end{theo}
\ni\ni{\it Proof.}\ \
(1)~It can be directly obtained from the proof of Proposition \ref{po1}.\\
(2)~Let $D_\lambda\in  CDer_{(\alpha,\beta,\gamma)}A$. Since $Cur \mathfrak{g}$ is a free $\mathcal{A}$-module of finite rank,
we can set $D_\lambda(a)=\sum_{i=0}^n\partial^id_\lambda^i(a)$ for every $a\in\mathfrak{g}$, where $d_\lambda^i$ are $\mathbb{C}$-linear maps from $\mathfrak{g}$ to $\mathfrak{g}[\lambda]$. According to (\ref{4000}),
we can get
\begin{equation}\label{err}
\alpha\sum_{i=0}^n\partial^id_\lambda^i([a,b])=\beta\sum_{i=0}^n(-\lambda-\mu)^i[d_\lambda^ia,b]
+\gamma\sum_{i=0}^n(\partial+\mu)^i[a,d_\lambda^i(b)].
\end{equation}
By Theorem \ref{r40003} and due to the simplicity of $A$, it is easy to see that
(i) and (ii) hold.

Next, let us consider the case when $\alpha=0$, $\beta=1$, $\gamma=-1$. By
(\ref{err}), we get $n=0$. Therefore, for any $D_\lambda \in CDer_{(0,1,-1)}A$,
we set $D_\lambda(a)=d_\lambda(a)$ for any $a\in \mathfrak{g}$, where $d_\lambda$ is a $\mathbb{C}$-linear map from $\mathfrak{g}$ to $\mathfrak{g}[\lambda]$. Thus, by (\ref{err}),
we get
\begin{equation}\label{err1}
[d_\lambda a,b]=[a, d_\lambda b].
\end{equation}
Set
$d_\lambda(a)=\sum_{i=0}^m\lambda^id_i(a)$
for every $a\in \mathfrak{g}$,
where $d_i$ are $\mathbb{C}$-linear maps from $\mathfrak{g}$ to $\mathfrak{g}$.
Taking it into (\ref{err1}), we obtain
\begin{equation}
[d_i(a),b]=[a,d_i(b)],~~\text{for all $i\in\{0,\cdots,m\}$}.
\end{equation}
It is known that $d_i=\alpha_i \text{Id}_{\mathfrak{g}}$, where $\alpha_i\in\C$ (see
Theorem 13 in \cite{ZZ1}). Therefore, $D_\lambda (a)=f(\lambda)a$ for any
$f(\lambda)\in \mathbb{C}[\lambda]$. Thus, (iii) holds.

Then, we consider the case when $\alpha=\delta$, $\beta=\gamma=1$.
Switching the place of $a$ and $b$ in (\ref{err}), we get
\begin{equation}\label{err2}
-\delta\sum_{i=0}^n\partial^id_\lambda^i([a,b])=\sum_{i=0}^n(-\lambda-\mu)^i[d_\lambda^ib,a]
+\sum_{i=0}^n(\partial+\mu)^i[b,d_\lambda^i(a)].
\end{equation}
Setting $\mu=0$ in (\ref{err}) and (\ref{err2}) and adding them up,
we obtain
\begin{equation}\label{err3}
\sum_{i=0}^n\partial^i([a,d_\lambda^i(b)]+[b,d_\lambda^i(a)])=\sum_{i=0}^n(-\lambda)^i([a,d_\lambda^i(b)]+[b,d_\lambda^i(a)]).
\end{equation}
Therefore, by (\ref{err3}), we get $[a,d_\lambda^i(b)]+[b,d_\lambda^i(a)]=0$ for any $i>1$ and $a$, $b\in \mathfrak{g}$. By Lemma 6.3 in \cite{DK},
$d_\lambda^i(a)=f_i(\lambda)a$, where $f_i(\lambda)\in \mathbb{C}[\lambda]$ for
$i>1$. Thus, $D_\lambda$ is of the following form:
\begin{equation}
D_\lambda(a)=d_\lambda^0(a)+P(\lambda,\partial)\partial a,
\end{equation}
where $P(\lambda,\partial)\in \mathbb{C}[\lambda,\partial]$.
Setting $P(\lambda,\partial)=\sum_{i=0}^mp_i(\lambda)\partial^i$ and by (\ref{err}), we have
\begin{eqnarray}
\delta(d_\lambda^0[a,b]+\sum_{i=0}^mp_i(\lambda)\partial^{i+1}[a,b])
&=&[d_\lambda^0(a),b]+\sum_{i=0}^m(-\lambda-\mu)^{i+1}p_i(\lambda)[a,b]\nonumber\\
\label{err4}&&+[a,d_\lambda^0(b)]+\sum_{i=0}^m(\partial+\mu)^{i+1}p_i(\lambda)[a,b].
\end{eqnarray}
Obviously, if $\delta=0$, by (\ref{err4}), $D_\lambda=0$. In the case of $\delta\neq 0$, if $m>1$, comparing the coefficients of $\partial^m$, we get
$(m+1)\mu p_m(\lambda)+(1-\delta)p_{m-1}(\lambda)=0$ which gives $p_m(\lambda)=0$. Therefore, we can set $D_\lambda(a)=d_\lambda^0(a)+p(\lambda)\partial a$. Taking it into (\ref{err4}), we get
\begin{equation}\label{err5}
\delta(d_\lambda^0[a,b])+\delta p(\lambda)\partial[a,b]
=[d_\lambda^0(a),b]+[a,d_\lambda^0(b)]-\lambda p(\lambda)[a,b]+\partial p(\lambda)[a,b].
\end{equation}

When $\delta=1$, $CDer_{(1,1,1)}A=CDer(A)=\{D_\lambda~|~D_\lambda(a)=P(\partial)(\partial+\lambda)a+d_\lambda(a),$
 where $a\in\mathfrak{g}$, $P(\partial)\in \mathbb{C}[\partial]$ and
    $d_\lambda\in CInn(A)\}$ (see Lemma 6.2 in \cite{DK}). When
    $\delta\neq 0$ and $\delta\neq 1$, by (\ref{err5}), we have $p(\lambda)=0$.
    Thus, in this case, $D_\lambda(a)=d_\lambda^0(a)$ and $d_\lambda^0(a)$ satisfies
  \begin{equation}\label{err6}
\delta(d_\lambda^0[a,b])
=[d_\lambda^0(a),b]+[a,d_\lambda^0(b)].
\end{equation}
It is known from Theorem 13 in \cite{ZZ1} that any $(\delta,1,1)$-derivation
of the Lie algebra $\mathfrak{g}$ is zero when $\delta\neq 1$ and $\delta\neq 2$, and any $(2,1,1)$-derivation of $\mathfrak{g}$ is equal to $\alpha \text{Id}_\mathfrak{g}$ for some complex number $\alpha$. Therefore,
by (\ref{err6}), we conclude that when $\delta\neq 1$ and $\delta\neq 2$,
$D_\lambda=0$ and when $\delta=2$, $D_\lambda(a)=f(\lambda)a$ for some
$f(\lambda)\in \mathbb{C}[\lambda]$. Thus, (iv) holds.

Finally, according to Theorem \ref{r40003}, $CDer_{(\delta,1,0)}A=CDer_{(0,1,-1)}A\cap CDer_{(2\delta,1,1)}A$. Then, (v)
can be easily obtained from (iii) and (iv).
\QED

\end{document}